\documentstyle[amsfonts,12pt]{article}
\input amssym.def

\newcommand{\h}{{\hbar}}
\newcommand{\bul}{{\bullet}}
\newcommand{\mb}{{\mathfrak{b}}}
\newcommand{\mg}{{\mathfrak{g}}}
\newcommand{\mh}{{\mathfrak{h}}}
\newcommand{\mL}{{\mathfrak{L}}}

\newcommand{\si}{{\sigma}}
\newcommand{\ga}{{\gamma}}

\newcommand{\pa}{{\partial}}

\newcommand{\cU}{{\cal U}}

\newcommand{\cO}{{\cal O}}

\newcommand{\cX}{{\cal X}}

\newcommand{\cK}{{\cal K}}

\newcommand{\bbK}{{\Bbb K}}
\newcommand{\bbC}{{\Bbb C}}
\newcommand{\bbR}{{\Bbb R}}

\newcommand{\bbQ}{{\Bbb Q}}
\newcommand{\bbF}{{\Bbb F}}

\newcommand{\te}{\theta}

\newcommand{\de}{{\delta}}
\newcommand{\D}{{\Delta}}
\newcommand{\ze}{{\zeta}}

\newcommand{\tH}{{\tilde{H}}}
\date{}
\newtheorem{defi}{Definition}
\newtheorem{lem}{Lemma}
\newtheorem{teo}{Theorem}

\newtheorem{prop}{Proposition}
\begin{document}
\begin{center}
{\huge\bf On The Formality Theorem for\\[0.3cm]
the Differential Graded \\[0.5cm] Lie Algebra  of Drinfeld}\\[0.3cm]
PoNing Chen\\[0.3cm]
{\it 77 Mass Ave. Cambridge, MA, 02139 \\[0.3cm]
Massachusetts Institute of Technology \\[0.5cm]
E-mail: pnchen@mit.edu}
\end{center}

\begin{abstract}
We discuss the differential graded Lie algebra (DGLA) of Drinfeld
modeled on the tensor algebra $\bigotimes \cU_{\mg}$ of the
universal enveloping algebra of a Lie algebra $\mg$ over any field
$\bbK$ of characteristic zero. We explicitly analyze
the first obstruction to the existence of
the formality quasi-isomorphism for this DGLA.
Our analysis implies the formality of
the DGLA $\bigotimes \cU_{\mb}$ of Drinfeld
associated to the two-dimensional Borel algebra $\mb$\,.
\end{abstract}

\section{Introduction}
Formality of a DG algebra (or equivalence of this DG
algebra to its homology in the homotopy category)
can be thought of as a ``no-go'' result which roughly
means that the DG algebra is poor as a homotopy
algebraic structure. On the other hand, formality
of a DG algebra provides us with remarkable
correspondences between homotopy invariant structures associated
to this DG algebra and its homology. These
correspondences give a lot of interesting applications
to deformation theory \cite{C}, \cite{D}, \cite{D1},
\cite{D2}, \cite{Gilles}, \cite{Kontsevich}, \cite{Sh} and \cite{TT}.

In this paper we discuss a DGLA which governs the
triangular deformations of Lie algebras
proposed and classified in \cite{D83} by
Drinfeld\footnote{See also
papers \cite{Peter}, \cite{Peter1}, \cite{Peter3},
and \cite{Peter2}, in which explicit triangular deformations
of various types of Lie algebras are discussed.}.
The total space of the
DGLA is the tensor algebra $\bigotimes \cU_{\mg}$ of the universal
enveloping algebra $\cU_{\mg}$ of a Lie algebra $\mg$
and the DGLA structure is defined by formulas
((\ref{D}),(\ref{d}),(\ref{brack})) in the next section.
With these definitions it is clear that triangular deformations
of a Lie algebra $\mg$ are in a bijective correspondence
with Maurer-Cartan elements of the DGLA $\bigotimes
\cU_{\mg}[[\h]]$\,.
Equivalent triangular deformations correspond to equivalent
Maurer-Cartan elements and vice-versa. For this reason
we refer to the DGLA structure ((\ref{D}),(\ref{d}),(\ref{brack})) on
$\bigotimes \cU_{\mg}$ as the DGLA of Drinfeld.

In paper \cite{C} by Calaque it is shown that for
any Lie algebra $\mg$ over $\bbR$ or $\bbC$ there
exists a quasi-isomorphism from the graded Lie algebra
$\bigwedge \mg$ to the DGLA $\bigotimes \cU_{\mg}$ of
Drinfeld associated to $\mg$.
This result, in particular, implies
Drinfeld's theorems about quantization of triangular
$r$-matrices \cite{D83}. However, Calaque's proof based on
Kontsevich's quasi-isomorphism \cite{Kontsevich} and
Dolgushev's globalization technique \cite{D} does
not apply to any field of characteristic zero.

In our paper we show
that there are no obstructions to the existence of the
first non-trivial structure map of an $L_{\infty}$
quasi-isomorphism from the $\bigwedge\mg$ to
$\bigotimes \cU_{\mg}$\,.
Our analysis implies the formality theorem for the DGLA
$\bigotimes \cU_{\mb}$ of Drinfeld associated to
the two-dimensional Borel algebra $\mb$\,.

The organization of the paper is as follows.
Section $2$ is preliminary. In this section we introduce
the notations and terminology we use in this paper.
In section $3$ we reduce the question of formality
of the DGLA of Drinfeld to an analysis of
admissible cocycles of the graded Lie algebra
$\bigwedge \mg$ with values in adjoint representation.
In the fourth section, we show that the $2$-cocycles
which appear as obstructions to the formality
of the DGLA of Drinfeld are exact. Using this fact, we
prove the formality theorem for the DGLA $\bigotimes\cU_{\mb}$
associated to the two-dimensional Borel
algebra $\mb\subset sl_2(\bbK)$\,. In the concluding section
we discuss recent results related to our
question.

\subsection*{Notation}

Throughout the paper, the underlying field $\bbK$ is assumed to be
of characteristic zero. We use $\Sigma$ to denote the suspension of
an graded vector space, namely, tensoring with $\bbK$ at degree $1$.
Given a cooperad ${\cO}$, we use $\bbF_{\cO}(V)$ to denote the
cofree coalgebra generated by $V$. Given a vector space $V$, we use
$Lie(V)$, $E(V)$ and $S(V)$ to denote the free Lie algebra, the
exterior algebra and the symmetric algebra generated by $V$,
respectively. Moreover, we use $S^n(V)$ to denote the subspace of
$S(V)$ consisting of  degree $n$ elements. We also use $V^{*}$ for
the dual vector space of $V$.

\subsection*{Acknowledgment}

I would like to express my sincere thanks to V. A. Dolgushev for
formulating the problem and for useful discussions of the topic.
This work was supported in part by
Department of Mathematics at MIT under the Summer
Program for Undergraduate Research in 2003 and by the Clay Mathematics Institute during summer 2005.

\section{Preliminaries.}
In this section, we introduce the notion of $L_{\infty}$-algebra
structures, the differential graded Lie algebra(DGLA) of Drinfeld
$\bigotimes \cU_{\mg}$, and give the motivation for our question.

Given a vector space $V$, we can construct $\bbF_{cocomm}(\Sigma^{-1}V)$ as the vector space
\begin{equation}
\bbF_{cocomm}(\Sigma^{-1}V)=\bigoplus_{n=1}^{\infty} S^{n}(\Sigma^{-1}V)\,
\end{equation}
with comultiplication $\D : \bbF_{cocomm}(\Sigma^{-1}V) \mapsto \otimes^2
\bbF_{cocomm}(\Sigma^{-1}V)$ given by
$$
\D(v_0)=0, \qquad \forall v_0 \in V
$$
and
$$
\D(v_0 \ldots v_n)=
\sum_{k=0}^{n-1} \sum_ {\te} sign(\te) \frac{1}{(n-k)!(k+1)!}v_{\te (0)} \ldots
v_{\te (k)} \otimes v_{\te(k+1)}  \ldots v_{\te(n)}
$$
where $\te$ is permutation of $n+1$ elements and $sign (\te)$ is the sign of permutation.
\begin{defi}
A coderivation of $ \bbF_{cocomm}(\Sigma^{-1}V)$ is a map
\begin{equation}
d: \bbF_{cocomm}(\Sigma^{-1}V) \mapsto  \bbF_{cocomm}(\Sigma^{-1}V)
\end{equation}
such that
\begin{equation}\label{Q}
\D d(X) =(d \otimes I+ I \otimes d) \D X, \qquad X \in \bbF_{cocomm}(\Sigma^{-1}V).
\end{equation}
We use $Coder( \bbF_{cocomm}(\Sigma^{-1}V))$ to denote the space of all coderivations of $\bbF_{cocomm}(\Sigma^{-1}V)$.
\end{defi}
We have the following proposition from \cite{GJ} that describes $Coder( \bbF_{cocomm}(\Sigma^{-1}V))$.
\begin{prop}(\cite{GJ} proposition 2.14)\label{GJ}
Let $\xi$ be the projection of  $\bbF_{cocomm}(\Sigma^{-1}V)$ to
$V$. The map sending $d\in Coder( \bbF_{cocomm}(\Sigma^{-1}V) $ to
$\xi \circ d \in Hom(\bbF_{cocomm}(\Sigma^{-1}V),V)$ is a vector
space isomorphism between $Coder( \bbF_{cocomm}(\Sigma^{-1}V))$ and
$Hom(\bbF_{cocomm}(\Sigma^{-1}V),V)$.
\end{prop}
As a result, a coderivation $d$ of degree $k$ is uniquely determined by a
semi-infinite collection of multi-linear maps
$$
d_n\,:\, S^{n} (\Sigma^{-1} V) \mapsto \Sigma^{-k-1} V\,
$$
where $d_n$ is the restriction of $\xi \circ d$ to $S^{n}
(\Sigma^{-1} V)$. For example,
$$
d(v)=d_1(v)\,, \qquad d(v_1, v_2)=d_2(v_1,v_2)+d_1(v_1)\otimes v_2 +(-1)^{(k_1-1)(k_2-1)}d_1(v_2)\otimes v_1
$$
$$
v \in V, \qquad v_1 \in V^{k_1}, \qquad v_2 \in V^{k_2}.
$$

Combining definition 4.2.14 \cite{GK} and proposition 4.2.15 \cite{GK}, we have the following definition:
\begin{defi}
A graded vector space $V$ is endowed with an $L_{\infty}$-algebra structure
if on $\bbF_{cocomm}(\Sigma^{-1}V)$, there is a coderivation $Q$ of degree
$1$ such that
\begin{equation}
\label{Q2}
 Q^{2}=0\,.
\end{equation}
\end{defi}
One can show that equation (\ref{Q2}) is equivalent to a collection of quadratic relations
on $Q_n$. The lowest two of these relations are
$$
Q_1^{2}=0\,,
$$
$$
Q_1(Q_2(v_1, v_2))=-Q_2(Q_1(v_1),v_2)+(-1)^{k_1}Q_2(v_1, Q_1(v_2)).
$$
{\bf Example.} Any differential graded Lie algebra (DGLA)
$(\mL,d,[,])$ is naturally an $L_{\infty}$-algebra with only two
non-vanishing maps
$$
Q_1 = d : \mL \mapsto \mL\,,
$$
$$
Q_2(h_1, h_2) = (-1)^{|h_1|} [h_1,h_2] : \otimes^2 \mL \mapsto \mL.
$$

\begin{defi}
An $L_{\infty}$-morphism $F$ between $L_{\infty}$-algebras $V$ and
$V^{\Box}$ is a homomorphism between  coalgebras
$\bbF_{cocomm}(\Sigma^{-1}V)$ and
$\bbF_{cocomm}(\Sigma^{-1}V^{\Box})$
\begin{equation}
\D F(X)=(F \otimes F)(\D X) \label{F},\qquad X\in \bbF_{cocomm}(\Sigma^{-1}V)
\end{equation}
such that
\begin{equation}
Q^{\Box} F(X) = F(Q X),\label{10} \qquad X \in \bbF_{cocomm}(\Sigma^{-1}V) .
\label{relation}
\end{equation}
\end{defi}

Similar to proposition \ref{GJ}, equation (\ref{F}) implies that $F$
is uniquely determined by its composition with projection onto
$V^{\Box}$. That is, $F$ is defined by a collection of maps
\begin{equation}
\label{maps} F_n : S^n (\Sigma^{-1} V) \mapsto \Sigma^{-1}V^{\Box}.
\end{equation}
For example,
\begin{equation}
F(v)=F_1(v)\,, \qquad F(v_1,v_2)=F_2(v_1, v_2) + F_1(v_1)\cdot  F_1(v_2).
\end{equation}

Equation (\ref{10}) can also be expressed in terms of relations on $F_n$ and $Q_n$.
These relations are more complicated for general $L_{\infty}$-algebras.
The first example is that $F_1$ commutes with the differentials,
\begin{equation}\label{motive}
F_1(Q_1 v)=Q_1^{\Box}F_1(v)\,.
\end{equation}
%This observation motivates the following natural definition.

{\bf Remark.}
For $L_{\infty}$-morphism between DGLAs, equation (\ref{10}) takes a simpler structure
since $Q_n$ is vanishing for all $n>2$. Let $F$ be an $L_{\infty}$-morphism  between DGLAs $(\mL, Q_1, Q_2)$
and $(\mL^{\Box}, Q^{\Box}_1, Q^{\Box}_2)$,
then equation (\ref{10}) takes the following form:
$$
Q^{\Box}_1 (F_n(\ga_1, \ga_2, \ldots, \ga_n))- \sum_{i=1}^n
(-1)^{k_1+\ldots+k_{i-1}+1-n} F_n(\ga_1, \ldots, Q_1 (\ga_i),
\ldots, \ga_n)
$$
\begin{equation}\label{MC}
=\frac12 \sum_{k,l\ge 1,~ k+l=n} \frac1{k!l!} \sum_{\si\in S_n}
\pm Q^{\Box}_2(F_k (\ga_{\si_1}, \ldots, \ga_{\si_k}), F_l (\ga_{\si_{k+1}},
\ldots, \ga_{\si_{k+l}})) \label{q-iso}
\end{equation}
$$
\sum_{i\neq j} \pm F_{n-1}(Q_2(\ga_i,\ga_j), \ga_1, \ldots,
\hat{\ga_i}, \ldots, \hat{\ga_j}, \ldots \ga_n), \qquad \ga_i \in
\mL^{k_i}_1\, .
$$
It is the form of the equation, rather than the actual signs in the equation, that is important to us.

Equation (\ref{motive}) motivates the following natural definition.
\begin{defi}
A quasi-isomorphism $F$ 
between $L_{\infty}$-algebras $V$ and $V^{\Box}$
is an $L_{\infty}$-morphism such that its first structure map, $F_1$, induces an isomorphism
between the cohomology spaces $H^{\bul}(V , Q_1)$ and
$H^{\bul}(V^{\Box} , Q^{\Box}_1)$.
\end{defi}

\begin{defi}
An element $r\in \mL^1 $ of a DGLA $\mL$ is a Maurer-Cartan element if it satisfies the following equation
\begin{equation}
dr+\frac{1}{2}[r,r]=0.
\end{equation}
\end{defi}

Quasi-isomorphisms between DGLAs are of primary importance in
deformation theory  because a quasi-isomorphism $F$ from $\mL$ to
$\mL^{\Box}$ gives a one-to-one correspondence between moduli spaces
of Maurer-Cartan elements of $\h\mL[[\h]]$ and $\h\mL^{\Box}[[\h]]$.
Given $F$, we define ${\tilde F}$ as follows:
\begin{equation}
\begin{array}{c}
{\tilde F}: \h \mL[[\h]] \mapsto \h \mL^{\Box}[[\h]]\\[0.3cm]
\displaystyle {\tilde F}( r)=\sum_{n=1}^{\infty}\frac{1}{n!} F_n(r,r,\ldots,r),  \label{twist}
\end{array}
\end{equation}
where $r$ is a Maurer-Cartan element of $\h \mL[[\h]]$.

It is easy to see that $\tilde F$ sends Maurer-Cartan elements to
Maurer-Cartan elements since $F$ is quasi-isomorphism. We refer the
readers to Chapter $2$ of  \cite{D2} for the proof that $\tilde F$
gives a bijection on moduli spaces of Maurer-Cartan elements.

\begin{defi}
A DGLA $\mL$ is formal if there is a quasi-isomorphism from its cohomology $H^{\bullet}(\mL)$ to $\mL$.
\end{defi}
\begin{defi}
We say that a graded Lie algebra $L$ is intrinsically formal if any DGLA $\mL$ with $H(\mL)=L$
is formal.
\end{defi}
As a result, formality theorem for a DGLA $\mL$ gives a one-to-one
correspondence between moduli spaces of Maurer-Cartan elements of
$\h\mL[[\h]]$ and $\h H^{\bullet}(\mL)[[\h]]$. Furthermore, in $\h
H^{\bullet}(\mL)[[\h]]$, equation for Maurer-Cartan elements is
simpler due to the  vanishing differential. An element $r\in \h
H^{\bullet}(\mL)[[\h]]$ is a Maurer-Cartan element if

\begin{equation}
[r,r]=0.
\end{equation}

In our paper, the main object of our discussion is the DGLA
$\bigotimes \cU_{\mg}$  whose DGLA structure controls the quantization of triangular $r$-matrices
proposed and described in paper \cite{D83} by Drinfeld.

\begin{defi}
A triangular $r$-matrix of a Lie algebra $\mg$ is an element $r\in
\wedge^2 \mg$ such that
\begin{equation}
[r,r]=0.
\end{equation}
\end{defi}
We are going to define the DGLA structure of $\bigotimes \cU_{\mg}$ as follows:

The symbol $\D$ denotes the standard comultiplication of
$\cU_{\mg}$, whose generators $I$ and $e\in \mg$ are acted on by
$\D$ as

\begin{equation}
\D I=I\otimes I, \qquad \D e= e \otimes I+ I
\otimes e. \label{comult}
\end{equation}
and $\D$ is a homomorphism $\cU_{\mg} \mapsto
\cU_{\mg}\otimes \cU_{\mg}$\,.

Using the comultiplication (\ref{comult}), we define the following
complex
\begin{equation}
\bigotimes \cU_\mg = D=\bigoplus_{k\ge {-1}} D^k, \qquad
D^{k}=\cU^{\otimes^{k+1}}_{\mg} \label{D}, \qquad D^{-1}=\bbK,
\end{equation}
with the differential $\de :D^k\mapsto D^{k+1}$ given by the
formula

\begin{equation}
\de \Phi = \Phi \otimes I - \sum_{l=0}^k(-)^{k-l} \D_{l}
\Phi+ (-)^k I \otimes \Phi, \qquad \Phi\in D^k, \label{d}
\end{equation}
where
$$
\D_l=I\otimes \ldots I\otimes
\underbrace{\D}_{the~l-th~place}\otimes I \ldots \otimes
I.
$$

In order to define a Lie algebra structure on the graded space
(\ref{D}), we define the following (non-associative) product

\begin{equation}
\Phi_1 \circ \Phi_2 = \sum_{i=0}^{k_1}(-)^{i k_2}
\tilde{\D}^{k_2}_i \Phi_1 (\Phi_2)_i, \qquad \Phi_1 \in
D^{k_1},\qquad \Phi_2 \in D^{k_2}, \label{prod}
\end{equation}
where
$$
(\Phi)_i=I\otimes \ldots I\otimes
\underbrace{\Phi}_{the~i-th~place} \otimes I \ldots \otimes I,
$$
$$
\D^k_i=I\otimes \ldots I \otimes
\underbrace{\D^k}_{the~i-th~place} \otimes I \ldots
\otimes I,
$$

and

$$
\D^k=(\underbrace{\D \otimes I \otimes \ldots
\otimes I}_{k}) (\underbrace{\D\otimes I\otimes \ldots
\otimes I}_{k-1}) \ldots (\D \otimes I) \D, \quad
k\ge 1,
$$
$$
\D^0=Id.
$$

A Lie bracket between homogeneous elements $\Phi_1 \in D^{k_1}$
and $\Phi_2 \in D^{k_2}$ is defined as

\begin{equation}
[\Phi_1,\Phi_2]_G=\Phi_1 \circ \Phi_2 - (-)^{k_1 k_2}\Phi_2 \circ
\Phi_1. \label{brack}
\end{equation}

The respective Jacobi identity can be proved by direct
computation.

It is easy to see that the differential (\ref{d}) is an inner
derivation of the Lie algebra (\ref{D}), (\ref{brack}),
\begin{equation}
\de \Phi= [I\otimes I, \Phi]_G \label{inder}.
\end{equation}
The Leibniz rule then follows from the Jacobi identity, and hence
the space (\ref{D}) is endowed with a DGLA structure.

We have the following PBW theorem:
\begin{teo}\label{PBW}
There exists a filtered vector space isomorphism $\sigma$ from
$\cU_{\mg}$ to $S(\mg)$.
\end{teo}
As a result, $\sigma$ induces a filtered vector space isomorphism
from $\bigotimes \cU_{\mg}$ to $\otimes S(\mg)$.

Notice that for a Lie algebra $\mg$, $E(\mg)$ has a
natural graded Lie algebra structure induced by extending the
Lie bracket of $\mg$ by the Leibniz rule with respect to exterior
product $\wedge$
\begin{equation}
[\ga, \ga_1\wedge \ga_2]=[\ga, \ga_1]\wedge \ga_2 +(-1)^{k
(k_1+1)}\ga_1\wedge [\ga,\ga_2]  \qquad \ga\in E^{k}(\mg)\,,~ \ga_1\in
E^{k_1}(\mg)\,. \label{extLie}
\end{equation}
We will use the same symbol for the Lie bracket in algebra $\mg$ and
the graded Lie bracket in algebra $E(\mg)$ defined above.

The structure of the cohomology space of the DGLA of Drinfeld is described
by the following theorem:

\begin{teo}\label{Drinfeldteo} \label{teo}
We have the following three statements:

(i) The cohomology space of the DGLA $D$ (\ref{D}), (\ref{d}), (\ref{brack}) is the
vector space of the exterior products of $\mg$, $E(\mg)$,
namely,
\begin{equation}
H^k(D)=\wedge^{k+1} \mg, \qquad k>-1; \qquad H^{-1}(D)=\bbK. \label{HKR}
\end{equation}

(ii) Let $f'$ be the map from $\otimes S(\mg)$ to $E(\mg)$ such that
\begin{equation} \label{ONO}
f'(g_1\otimes \dots \otimes g_k)= g_1\wedge \dots \wedge g_k
\end{equation}
where $g_i \in \mg$ and that $f(a_1\otimes\dots\otimes a_s)=0$ where $a_j \in S(\mg)$ if for some
$i$, $a_i \in S(\mg)$ has degree other than $1$. Then,
\begin{equation}\label{isom1}
f= f' \circ \sigma
\end{equation}
is a quasi-isomorphism from $\bigotimes \cU_{\mg}$ to $E(\mg)$
where $\sigma$ is the vector space isomorphism in theorem \ref{PBW}.

(iii) The Lie algebra structure of $E(\mg)$ induced by that of $D$ is the same as
the one extended from bracket of $\mg$ by equation (\ref{extLie}).
\end{teo}
{\bf Proof.} For the first statement, we may assume that $\mg$ is
commutative since the definition of differential of $D$ uses only
the comultiplication  of $\cU_{\mg}$ but not the multiplication of
$\cU_{\mg}$. Thus, $\cU_{\mg}=S(\mg)$. We construct an isomorphism
$i$ from  the complex $(\otimes S(\mg), \D)$ to the Hochschild
cochain complex $C^{\bullet}(S(\mg^{*}),\bbK)$. Pick a basis
$\{g_1,\dots g_n\}$ for $\mg$. We have
\begin{equation}
S(\mg)=\bbK[g_1,\dots ,g_n]
\end{equation}
and
\begin{equation}
S(\mg^{*})=\bbK[g_1^{*},\dots,g_n^{*}].
\end{equation}
Given $b\in \bbK[g_1^{*},\dots,g_n^{*}]$, we set
\begin{equation}
i(g_{k_1}\dots g_{k_s})(b)=\frac{\pa}{\pa g_{k_1}^{*}}\dots \frac{\pa}{\pa g_{k_s}^{*}}(b)|_{g*=0}.
\end{equation}
This defines a map $i: S(\mg) \mapsto Hom(S(\mg^{*}),\bbK)$ and we extend this map to
\begin{equation}
i :\otimes S(\mg) \mapsto C^{\bullet}(S(\mg^*),\bbK)
\end{equation}
by the formula
\begin{equation}
i(a_1 \otimes\dots\otimes a_k)(b_1\otimes\dots\otimes b_k)=i(a_1)(b_1)\times\dots\times i(a_k)(b_k)
\end{equation}
where $a_i \in S(\mg)$ and $b_i \in S(\mg^*)$.

In order to compute $HH^{\bullet}(S(\mg^*),\bbK)$, we use the
following equation
\begin{equation}
HH^{\bullet}(A,\bbK)=Ext^{\bullet}_{A\otimes A^{op}}(A,\bbK)
\end{equation}
where $A$ is an augmented algebra and acts on $\bbK$ by multiplying the image of the augmentation.

Let $A=S(\mg^*)=\bbK [x_1,\dots,x_n ]=\bbK [X]$. Then $A\otimes A=\bbK [X]\otimes \bbK [Y]$ and let $\theta^i$ be
anti-commutative variables $1 \le i \le n$. We use the Koszul
resolution $\cK^{\bullet}$ of $A$ to compute
$Ext^{\bullet}_{A\otimes A}(A,\bbK)$.
\begin{equation}
\cK^{\bullet}=\bbK [ X, Y,  \theta]
\end{equation}
where $\cK^{i}$ is the module consisting of elements of degree $i$
in $\theta$ and with differential
\begin{equation}
\pa_{\cK}=\sum_{j}(x^j-y^j)\frac{\pa}{\pa \theta^{j}}.
\end{equation}

Direct computation shows that
\begin{equation}
Hom_{A\otimes A} (\cK^m, \bbK)= Hom_{\bbK} (E^m(\mg^{*}),\bbK)=E^m(\mg)
\end{equation}
and that differential on $Hom_{A\otimes A} (\cK^{\bullet}, \bbK)$ induced by $\pa_K$ is zero.
This proves the first statement of the theorem.

For the second statement, we use the following quasi-isomorphism
between Koszul resolution and bar resolution of $A$:
\begin{equation}\label{qkb}
c_{i_1\dots i_k} a \otimes \theta^{i_1}\dots \theta^{i_k}\otimes b \mapsto
c_{i_1\dots i_k} a \otimes x^{i_1} \otimes \dots \otimes x^{i_k}\otimes b
\end{equation}
where $a\otimes b \in A\otimes A$. Notice that $\theta$s are anti-commutative variable
and $x,y$ s are commutative variable. Let $f$ be the map defined in equation (\ref{isom1}),
and $h$
\begin{equation}\label{isom2}
h : E(\mg) \mapsto \bigotimes \cU_{\mg}
\end{equation}
be the natural map from $E(\mg)$ to $D$. It is easy to see that $f \circ h$
is the identity map. Furthermore, using the quasi-isomorphism (\ref{qkb}), we can also show that
$h \circ f$ is homotopic to the identity. Thus we can conclude that $f$ is a quasi-isomorphism.
Similarly, we can proves (iii) using the maps $f$ and $h$. $\qquad \Box$

Since $H(\bigotimes \cU_{\mg})=E(\mg)$ and a triangular $r$-matrix
$r\in \wedge^2 (\mg)= H^1(\bigotimes \cU_{\mg})$ is also a
Maurer-Cartan element, a quasi-isomorphism $F$ from $E(\mg)$ to
$\bigotimes \cU_{\mg}$ induces a bijection of the moduli spaces of
Maurer-Cartan elements of $\h E(\mg)$ and $\h (\bigotimes
\cU_{\mg})$ given by equation (\ref{twist}).  Moreover, we can
construct

\begin{equation}\label{deform1}
T=I\otimes I +\tilde F(r).
\end{equation}
Since $\tilde F(r)$ is a Maurer-Cartan element, $T$ satisfies the following equation:

\begin{equation}\label{CYB}
T_{12,3}T_{12}=T_{1,23}T_{23}.
\end{equation}

The construction of an element $T$ satisfying equation (\ref{CYB})
is called the quantization of triangular $r$-matrices \cite{D83}.

\section{Obstructions to the formality of $L_{\infty}$-Algebras}
Given an $L_{\infty}$-algebra $(\mL,Q)$, its cohomology $H=H(\mL)$
is naturally a DGLA with $Q^H_1=0$ and $Q^H_2([a],[b])=[Q^H(a,b)]$
where  $a,b \in \mL$ represent cohomology classes $[a]$ and $[b]$ in
$H(\mL)$. Denote this DGLA by  $(H,Q^H)$. Let $(\mh,d)$ the DGLA of
coderivation of $\bbF_{cocomm}(\Sigma^{-1} H)$ with differential defined by commutator
with $Q^H$. In this section we show that obstruction to the
formality of $\mL$ lies in $H^1(\mh)$. A more general result is
proven in section $4$ of \cite{Hi}. However, we decide to give a
proof without referring to the technique of closed model category
\cite{Q} and to give an explicit construction.
\begin{lem}\label{4.2.1}
Given a quasi-isomorphism of complex $f$ form $H$ to $\mL$, there
exists an $L_{\infty}$ structure $\tH$ on the space $H$ such that
$\tH$ is quasi-isomorphic to $\mL$ and that the first structure map
of the quasi-isomorphism is $f$.
\end{lem}
{\bf Proof.}
We will prove the lemma by constructing a quasi-isomorphism $F$ from
$(H[\lambda],\tilde Q)$ to $(\mL[\lambda],\tilde Q^H)$ where
${\tilde Q}$ and ${\tilde Q^H}$ on $\mL[\lambda]$ and $H[\lambda]$ are
$$
{\tilde Q_n} = \lambda^{n-1} Q_n,
$$
and
$$
{\tilde Q^{H}_n} =  \lambda^{n-1} Q'_n
$$
and $F$ have the following form:
$$
{\tilde F_n} = \lambda^{n-1} F_n
$$
and
$$
{\tilde F_1} = f.
$$
We construct ${\tilde Q^{H}_n}$ and ${\tilde F_n} $ inductively as follows:
For $n=2$, notice that $f=F_1$ is quasi-isomorphism from $H$ to $\mL$. As a result,
\begin{equation}
F_1(Q_2^H(\ga_1,\ga_2))-Q_2(F_1(\ga_1),F_1(\ga_2))
\end{equation}
is $Q_1$ exact. Hence, there exists $F_2$ such that
\begin{equation}
F_1(Q_2^H(\ga_1,\ga_2))-Q_2(F_1(\ga_1),F_1(\ga_2))=Q_1(F_2(\ga_1,\ga_2)).
\end{equation}
Suppose that we have defined ${\tilde Q_k^{H}}$ and ${\tilde F_k}$ such that
\begin{equation}\label{ind}
{\tilde Q}{\tilde F} - {\tilde F} {\tilde Q^H} = 0 \,\,\,mod\, \lambda^{n-1}
\end{equation}
and that
\begin{equation}\label{vanish}
{\tilde Q^H}{\tilde Q^H} = 0\,\,\, mod \,\lambda^{n},
\end{equation}
where $n > 3$. Notice that equation (\ref{ind}) implies that
\begin{equation}\label{n+1}
{\tilde Q}({\tilde Q}{\tilde F} - {\tilde F} {\tilde Q^H} )+
({\tilde Q}{\tilde F} - {\tilde F} {\tilde Q^H}) {\tilde Q^H}=0 \,\,\, mod\, \lambda^{n}.
\end{equation}
Let $\Phi_n$ be the sum of terms of degree $(n-1)$ in $\lambda$ of
$\xi({\tilde Q}{\tilde F} - {\tilde F} {\tilde Q^H} )$
where $\xi$ is the projection of $\bbF_{cocomm}(\Sigma^{-1}H)$ to $H$.
Equation (\ref{n+1}) implies that
\begin{equation}
Q_1 (\Phi_n)=0.
\end{equation}
As a result, there exists ${\tilde F_n}$ and ${\tilde Q^H_n}$ such that
\begin{equation}\label{step}
\Phi_n + Q_1(\tilde F_n)-F_1(\tilde Q^H_n)=0.
\end{equation}
Let's denote the change of modifying $\tilde Q^H$ and $\tilde F$ with $\tilde Q^H_n$ and
$\tilde F_n$ by $\delta \tilde Q^H$ and $\delta \tilde F$.
We have
\begin{equation}
({\tilde Q}{\tilde F} - {\tilde F} {\tilde Q^H})+\tilde Q \delta \tilde F -
\delta \tilde F \tilde Q^H -\tilde F \delta \tilde Q^H=0 \,\,\,  mod\, \lambda^{n}.
\end{equation}
Furthermore, we have that
$$
{\tilde F} (\delta \tilde Q^H \tilde Q^H + \tilde Q^H+\delta \tilde Q^H)
$$
$$
={\tilde Q} {\tilde F}{ \tilde Q^H} -{ \tilde F}
({\tilde
Q^H})^2+{\tilde Q} \delta {\tilde F}{ \tilde Q^H} + {\tilde Q}{
\tilde F} \delta {\tilde Q^H} \,\,\,mod\, \lambda ^{n+1}
$$
$$
=- \tilde F (\tilde Q^H)^2 + \tilde Q (\lambda ^{n}g) \,\,\, mod\, \lambda ^{n+1}
$$
where $g$ is of degree $0$ or higher in $\lambda$ which shows that

$$
\xi \tilde Q (\lambda ^{n}g) = Q_1 ( (\lambda ^{n}g) \,\,\, mod\, \lambda ^{n+1}.
$$
This shows that
$$
F_1 (\xi((\delta \tilde Q^H \tilde Q^H + \tilde Q^H+\delta \tilde Q^H + (\tilde Q^H)^2)) ) +
Q_1 ( (\lambda ^{n}g) =0 \,\,\, mod\, \lambda^{n+1}.
$$
Hence
$$
 \delta \tilde Q^H \tilde Q^H + \tilde Q^H+\delta \tilde Q^H + \tilde (Q^H)^2)=0 \,\,\, mod\, \lambda^{n+1}.
$$
That is,
$$
\xi (\tilde Q^H+\delta \tilde Q^H)^2=0 \,\,\, mod\, \lambda^{n+1},
$$
since $(\delta \tilde Q^H)^2$ is of degree $2n-2 \ge n+1$ in $\lambda$. This finishes the proof. $\qquad \Box$

Now we want to construct a quasi-isomorphism from $H$ to $\tilde H$.
Hinich uses the following inductive procedure to this construction.
Let $H[\h]$ be the DGLA with
\begin{equation}
Q^{\h}_n =\h ^{n-1} Q^H_n.
\end{equation}
Suppose we have constructed  a map $g: \bbF_{cocomm}(\Sigma^{-1} H[\h]) \mapsto \bbF_{cocomm} (\Sigma^{-1} \tilde H[\h])$ whose structure maps are
\begin{equation}\label{degree}
g_n= \h^{n-1} \phi_n
\end{equation}
where $\phi_n$ are structure maps for $L_{\infty}$ morphism from $ H$ to  $ \tilde H$
such that
\begin{equation}\label{pro}
 g^{*} : \bbF_{cocomm}(\Sigma^{-1}(H[\h]/(\h^{n}))) \mapsto \bbF_{cocomm}(\Sigma^{-1} ( \tilde H[\h]/(\h^{n})))
\end{equation}
is a quasi-isomorphism. Due to lemma \ref{4.2.1}, there exists a codifferential $Q'$ on \\
$\bbF_{cocomm}(\Sigma^{-1}(H[\h]/(\h^{n+1})))$ such that
\begin{equation}
 g^{*} : (\bbF_{cocomm}(\Sigma^{-1} (H[\h]/(\h^{n+1}))),Q') \mapsto \tilde \bbF_{cocomm}(\Sigma^{-1} (H[\h]/(\h^{n+1})))
\end{equation}
is a quasi-isomorphism. In fact, $Q^{\h}_k=Q'_k$ for $k \le n$ due
to equation (\ref{pro}). Thus we conclude that $Q'_{n+1}$ viewed as
a coderivation is a closed element of $H^1(\mh)$ since $(Q')^2=0$.
If $Q'_{n+1}$ is trivial in $H^1$, say $Q'_{n+1}= d z$, then we let
$j$ be the $L_{\infty}$ morphism with two non-vanishing structure
maps $j_1=Id$ and $j_{n+1}=z$ then the maps $g'=g\circ j$ satisfies
equation (\ref{degree}) and induces quasi-isomorphism from
$H[\h]/(\h^{n+1})$ to $\tilde H[\h]/(\h^{n+1})$.

Combining this procedure with lemma \ref{4.2.1}, we prove the particular case of Theorem 4.2 in \cite{Hi}.

\begin{teo}
Given an $L_{\infty}$-algebra $(\mL,Q)$ with cohomology $(H,Q^{H})$, let $\mh$ be the DGLA
of coderivation of coalgebra $\bbF_{cocomm}(\Sigma^{-1}H)$ with differential defined by
commutator with $Q^{H}$.  If $H^1(\mh)=0$ then $H$ is intrinsically formal.
\label{Hinich}
\end{teo}

\section{Vanishing of the first obstruction to formality of $\bigotimes \cU_{\mg}$}
In this section, we use the proof of theorem \ref{Hinich} to show that the first obstruction to
formality of $\bigotimes \cU_{\mg}$ vanishes.

Direct computation shows that $H^1(\mh)$ is nothing but the Lie algebra cohomology with the following differential
\begin{equation}\label{48}
\pa \Phi (\ga_0,\ga_1 \ldots \ga_{n+1})= - \sum_{i=0}^{n+1}
(-1)^{i+k_i(k_0+ \ldots k_{i-1})}[\ga_i,\Phi(\ga_0,\ga_1 \ldots
\tilde{\ga_i} \ldots \ga_{n+1})]
\end{equation}
$$
- \sum_{i=0}^{n}\sum_{j=i+1}^{n+1} (-1)^{i+j+k_i(k_0+ \ldots
k_{i-1})+k_j(k_0+ \ldots k_{j-1}-k_i)}\Phi([\ga_i,\ga_j],\ga_0
\ldots \tilde{\ga_i} \ldots \tilde{\ga_j} \ldots \ga_{n+1})\,,
$$
where $\ga_i$ has degree $k_i$ (Example 4.2.4 b \cite{GK}).

\begin{defi}
A cocycle $\Phi: \wedge^{n+1} E(\mg)\mapsto E(\mg)$ is admissible if any component of $\Phi(\ga_0,
\ldots  ,\ga_n)$ is expressed in terms of components of $\ga_i$
via the Lie bracket of $\mg$ and the coefficients depend only on
the degrees of poly-vectors $\ga_i$ and that the cochain is a
cocycle for any Lie algebra $\mg$.
\end{defi}

We have the following proposition for the cochains showing up as obstruction to the formality of $\bigotimes \cU_{\mg}$.
\begin{prop}
If one uses the pair of quasi-isomorphisms $f$ and $h$ defined in equation (\ref{isom1}) and equation(\ref{isom2})
then obstructions to the formality of $\bigotimes \cU_{\mg}$ are represented by admissible $n$-cocycles of degree $2-n$.
\end{prop}
{\bf Proof.}
Let $h$ and $f$ be the quasi-isomorphisms defined in equation (\ref{isom1}) and equation (\ref{isom2}).
We follow the construction in lemma \ref{4.2.1} with $F_1=h$.
From equation (\ref{step}) we know that the first non-vanishing $Q_n(\ga_1,
\ldots, \ga_n)$ is a cocycle of $E(\mg)$ representing the cohomological class of
\begin{equation}\label{49}
\begin{array}{c}
\displaystyle \frac12 \sum_{k,l\ge 1,~ k+l=n} \frac1{k!l!}
\sum_{\si\in S_n} \pm [F_k (\ga_{\si_1}, \ldots, \ga_{\si_k}), F_l
(\ga_{\si_{k+1}}, \ldots, \ga_{\si_{k+l}})]_G \label{q-isoQQQ}
\\[0.3cm]
\displaystyle -\sum_{i\neq j} \pm F_{n-1}([\ga_i,\ga_j], \ga_1,
\ldots, \hat{\ga_i}, \ldots, \hat{\ga_j}, \ldots \ga_n).
\end{array}
\end{equation}
Let the above sum be $\ze$. Consider the element $f(\ze)$. It is a cocycle in $E(\mg)$
and $h(f(\ze))$ represents the cohomology class of $\ze$. Thus, in equation (\ref{step}) we may choose
$Q_n$ to be $f(\ze)$. Furthermore, from the construction in lemma \ref{4.2.1}, we see that $F_k$ can be chosen
so that any component of the $f(\ze)$ resulting from the construction is expressed in terms
of components of $\ga_i$ via the Lie bracket of $\mg$ and that the
coefficient is independent of the particular choice component of
$\ga_i$. Thus, $Q_n$ is a admissible $(n-1)$-cocycle of degree $2
- n$.  $\qquad \Box$

From the above proposition and theorem \ref{Hinich}, we get:
\begin{teo}
If all obstructions to formality of $\bigotimes \cU_{Lie(V)}$
are trivial for any finite dimensional vector space $V$, then for
any finite dimensional Lie algebra $\mg$, $\bigotimes \cU_{\mg}$ is formal.
\end{teo}

In the rest of this section,  we prove that the first obstruction for the formality of $\bigotimes \cU_{Lie(V)}$ vanishes.
Thus we prove that  $\bigotimes \cU_{\mb}$ is formal where $\mb$ is a two-dimensional  Lie algebra.

\begin{teo}
The first obstructions to formality of $\bigotimes \cU_{Lie(V)}$ are trivial.
\end{teo}
{\bf Proof.} First we recall that first obstructions are admissible
$2$-cocycles of degree $-1$. Let the first obstruction to be
\begin{equation}
\Phi (\gamma_0,\gamma_1,\gamma_2)
\end{equation}
which is a summation of terms which are expressed as wedge products of components of the polyvectors $\gamma_i$ via
three bracket operations. Three bracket operations can involve four to six components of $\gamma_i$. We start from
analyzing terms with the least number of components involved.

From now on, we will assume $\gamma_i$ to be polyvectors given by
\begin{equation}
 \gamma_0=x_0 \wedge x_1 \wedge \dots \wedge x_p,
\end{equation}
\begin{equation}
 \gamma_1=y_0 \wedge y_1 \wedge \dots \wedge y_q,
\end{equation}
\begin{equation}
 \gamma_2=z_0 \wedge z_1 \wedge \dots \wedge z_r,
\end{equation}
\begin{equation}
 \gamma_3=w_0 \wedge w_1 \wedge \dots \wedge w_s.
\end{equation}
Up to permutation of polyvectors and the indices, there are four types of terms involving four components of $\gamma_i$.
They are\\[0.5cm]
\begin{equation}\label{type1}
[[[x_0,y_0],x_1],x_2]\wedge x_3\dots \wedge x_p \wedge y_1 \dots \wedge y_q \wedge \gamma_2
\end{equation}
\begin{equation}\label{type2}
[[[x_0,y_0],x_1],y_1]\wedge x_2\dots \wedge x_p \wedge y_2 \dots \wedge y_q \wedge \gamma_2
\end{equation}
\begin{equation}\label{type3}
[[[x_0,y_0],x_1],z_0]\wedge x_2\dots \wedge x_p \wedge y_1 \dots \wedge y_q \wedge z_1\dots \wedge z_r
\end{equation}
\begin{equation}\label{type4}
[[[x_0,y_0],z_0],x_1]\wedge x_2\dots \wedge x_p \wedge y_1 \dots \wedge y_q \wedge z_1\dots \wedge z_r
\end{equation}
We will show that in the cocycle, terms of the types (\ref{type1})
and (\ref{type2}) do not appear and terms of types (\ref{type3}) and
(\ref{type4}) can be killed by adding coboundary terms.

First we consider terms in $\partial \Phi$ with brackets involving
$x_0$ , $x_1$ , $x_2$ , $y_0$ , $w_0$. There are four sums
containing such terms. Two of them are coming from the sums
$[\Phi(\gamma_0,\gamma_1,\gamma_2),\gamma_3]$ and
$[\Phi(\gamma_0,\gamma_2,\gamma_3),\gamma_1]$ and two of them are
from the sums $\Phi([\gamma_0,\gamma_1],\gamma_2,\gamma_3)$ and
$\Phi([\gamma_0,\gamma_3],\gamma_2,\gamma_1)$. Furthermore, we can
notice that the first two are antisymmetric when we exchange $x_0$
with $x_1$ but the last two are symmetric. Thus in a cocycle, the
contribution of terms involving $x_0$, $x_1$, $x_2$, $y_0$, $w_0$
from the first two sums must be zero. That is we have
\begin{equation}
\alpha_1 [[[[x_0,y_0],x_1],x_2],w_0]+\alpha_2 [[[[x_0,w_0],x_1],x_2],y_0]=0,
\end{equation}
where $\alpha_1$ and $\alpha_2$ are coefficients of the terms in the first two sums respectively.
However, $\mathfrak{g}$ is free Lie algebra and thus
\begin{equation}
\alpha_1 =\alpha_2 =0,
\end{equation}
which shows that terms of the type  (\ref{type1}) do not appear.

Secondly, we consider terms in $\partial \Phi$ with brackets
involving $x_0$ , $x_1$ , $y_0$ , $y_1$ and $w_0$ which are
antisymmetric in the pairs ($x_0$,$x_1$) and ($y_0$,$y_1$). In,
$\partial \Phi$, there is only one such term from the sum
$[\Phi(\gamma_0,\gamma_1,\gamma_2),\gamma_3]$ and thus in a cocycle
$\Phi(\gamma_0,\gamma_1,\gamma_2)$, terms of  type (\ref{type2}) do
not appear.

Finally, we will show that terms of types (\ref{type3}) and (\ref{type4}) can be killed by
adding coboundary terms $\partial \phi(\gamma_0,\gamma_1,\gamma_2)$ where
\begin{equation}
\phi(\gamma_0,\gamma_1)=\sum a_{i,j,k} [[x_i,y_j],x_k]\wedge \tilde{\gamma_0} \wedge \tilde{\gamma_1}
\end{equation}
in which $\tilde{\gamma_0}$ and $\tilde{\gamma_1}$ are obtained from the original
polyvectors by removing the components appeared in front.

We can choose $a_{i,j,k}$ so that the terms of the type
(\ref{type3}) are killed by adding $\partial
\phi(\gamma_0,\gamma_1,\gamma_2)$. That is, in $(\Phi-\partial
\phi)(\gamma_0,\gamma_1,\gamma_2)$, the only possible terms
involving four components in the brackets are the fourth type. Now,
we consider terms in $\partial \Phi$ with brackets involving $x_0$ ,
$x_1$ , $y_0$ , $z_0$ and $w_0$ and are antisymmetric when
permuting $x_0$ with $x_1$.  There are three such terms from the sums
$[\Phi(\gamma_0,\gamma_1,\gamma_2),\gamma_3]$ and cyclic permutation of $1,2,3$.\\
That $\Phi$ is a cocycle  implies
\begin{equation}
\alpha_{y,z}[[[[x_0,y_0],x_1],z_0],w_0]+\alpha_{z,w}[[[[x_0,z_0],x_1],w_0],y_0]+
\alpha_{y,w}[[[[x_0,y_0],x_1],w_0],z_0]=0.
\end{equation}
Again, $\mathfrak{g}$ is free Lie algebra and thus all three $\alpha$'s above are zero.
As a result, up to a coboundary, terms with four components in brackets are zero.

Next we consider terms with five components in brackets. Up to permutation of polyvectors
and the indices, there are two types of terms\\[0.5cm]
\begin{equation}\label{type5}
[[x_0,y_0],z_0] \wedge [*,*']\wedge \tilde{\gamma_0} \wedge \tilde{\gamma_1} \wedge \tilde{\gamma_2}
\end{equation}
\begin{equation}\label{type6}
[[x_i,y_j],x_k] \wedge [*,*']\wedge \tilde{\gamma_0} \wedge \tilde{\gamma_1} \wedge \tilde{\gamma_2}
\end{equation}
where $\tilde{\gamma_i}$ is obtained from original polyvectors by removing the components appeared in front.

First we consider terms in $\partial \Phi$ of the form:
\begin{equation}
[[[x_0,y_0],z_0],w_0] \wedge [*,*']\wedge \tilde{\gamma_0} \wedge \tilde{\gamma_1} \wedge \tilde{\gamma_2} \wedge \tilde{\gamma_3}.
\end{equation}
In general, we may assume $*=x_1$, there are at most three terms in $\partial \Phi$ of the above form
and is symmetric when exchanging $x_0$ and $x_1$ which are coming from the sum
$\Phi([\gamma_0,\gamma_1],\gamma_2,\gamma_3)$ and cyclic permutation of (1,2,3). Thus we have
\begin{equation}
\beta_{x,y,z}[[[x_0,y_0],z_0],w_0]+\beta_{x,z,w}[[[x_0,z_0],w_0],y_0]+\beta_{x,y,w}[[[x_0,y_0],w_0],z_0]=0,
\end{equation}
which implies $\beta$ are all zero since $\mathfrak{g}$ is free Lie
algebra. As a result, the terms of type (\ref{type5}) do not appear
in $\Phi-\partial \phi$.

For the terms of type (\ref{type6}), if one of the $*$ or $*'$ is from the polyvector $\gamma_0$ or $\gamma_1$ then
we can use exactly the same argument as above to show that such terms do not appear in cocycles. Thus we only need to
consider the term
\begin{equation}
[[x_0,y_0],x_1] \wedge [z_0,z_1]\wedge \tilde{\gamma_0} \wedge \tilde{\gamma_1} \wedge \tilde{\gamma_2}.
\end{equation}
In $\partial \Phi $ we consider terms of the form
\begin{equation}
[[[x_0,y_0],x_1],w_0] \wedge [z_0,z_1]\wedge \tilde{\gamma_0} \wedge \tilde{\gamma_1} \wedge \tilde{\gamma_2}\wedge \tilde{\gamma_3},
\end{equation}
which are symmetric with respect to $x_0$ and $x_1$. There are at most two such terms from the sums
$\Phi([\gamma_0,\gamma_1],\gamma_2,\gamma_3)$ and $\Phi([\gamma_0,\gamma_3],\gamma_2,\gamma_1)$. And for a coboundary, we have
\begin{equation}
a[[[x_0,y_0],w_0],x_1]+b[[[x_0,w_0],x_0],x_1]=0,
\end{equation}
which implies $a$ and $b$ are zero and we finish the proof of the statement that up to a coboundary, terms
with five components in brackets do no appear in the cocycle.

Finally we consider terms with six components in brackets. Up to permutation of
polyvectors and the indices, there are four possible types of terms\\[0.5cm]
\begin{equation}\label{type7}
[x_0,y_0] \wedge [x_1,z_0] \wedge [y_1,z_1] \wedge \gamma_0 \wedge \gamma_1 \wedge \gamma_2
\end{equation}
\begin{equation}\label{type8}
[x_0,x_1] \wedge [x_2,y_0] \wedge [y_0,z_0] \wedge \gamma_0 \wedge \gamma_1 \wedge \gamma_2
\end{equation}
\begin{equation}\label{type9}
[x_0,x_1] \wedge [x_2,y_0] \wedge [x_3,z_0] \wedge \gamma_0 \wedge \gamma_1 \wedge \gamma_2
\end{equation}
\begin{equation}\label{type10}
[x_0,x_1] \wedge [y_0,y_1] \wedge [z_0,z_1] \wedge \gamma_0 \wedge \gamma_1 \wedge \gamma_2
\end{equation}
Notice that in the first three of  the above terms, there exists one bracket such that for each component in the
bracket, the polyvector of that components have another components in another bracket.
$[x_1,z_0]$, $[x_2,y_0]$ and $[x_0,x_1]$ are examples of such brackets in the three terms respectively.

Let's first consider terms of the form
\begin{equation}
[x_0,y_0], \wedge [[x_1,z_0],w_0] \wedge \tilde{\gamma_0} \wedge \tilde{\gamma_1} \wedge \tilde{\gamma_2}\wedge \tilde{\gamma_3}
\end{equation}
in $\partial \Phi$ which is anti-symmetric with respect to the pair
$(x_0,x_1)$ and $(z_0,z_1)$. There is only one such term from the
sum $[\Phi(\gamma_0,\gamma_1,\gamma_2),\gamma_3]$ and thus for a
cocycle, the terms of type (\ref{type7}) are vanishing.

Using exactly the same argument, we can show that terms of type (\ref{type8}) and (\ref{type9}) are vanishing.

For terms of type (\ref{type10}), it is clear from equation
(\ref{49}) that such terms does not show up as obstructions in our
construction using the pair of quasi-isomorphism $f$ and $h$
((\ref{isom1}),(\ref{isom2})). Thus, we have shown that the first
obstruction is trivial. $\qquad \Box$

As we mentioned in previous section, if we can prove all
obstructions are trivial, then we prove formality theorem. However,
we may notice that the result for every two-cocycle implies
the formality theorem for two dimensional Lie algebra since we have
the following:
\begin{lem}
For a two-dimensional Lie algebra $\mg$, all admissible
$n$-cochains of degree $1-n$ for $n>2$ are vanishing.
\end{lem}
{\bf Proof.} We know that any non-abelian two-dimensional Lie
algebra, up to isomorphism, is $\mb$, the Borel subalgebra of
$sl_2(\bbK)$, namely
\begin{equation}
[e_1,e_2]=2e_1\,.
\end{equation}
Also notice that the space $E(\mg)$ is spanned by $e_1$,$e_2$,$e_1
\wedge e_2$ and thus, for an $n$-cochain, we have at most $n+1$
components of $e_2$. However, we may notice that for each bracket,
we must have a component of $e_2$ or the bracket will otherwise be
zero. For a non-zero $n$-cochain of degree $1-n$, we have $2n-1$
bracket. As a result, or $n>2$ the number of brackets is more than
the possible number of components of $e_2$. Therefore, we prove
that all $n$-cochains $n>2$ are vanishing. $\qquad \Box$

Thus we have proved the following:
\begin{teo}
For the two-dimensional Borel Lie algebra $\mb$,
$\bigotimes \cU_\mb$ is formal.
\end{teo}

\section{Concluding remarks}
We would like to mention that the construction
of the cohomology classes of $E(\mg)$
which appear as obstructions to the
formality depends on the choice of the
quasi-isomorphism (\ref{isom1}).
Therefore, these classes are not homotopy
invariant. However, the classes can be used to introduce
a partially defined products, so-called
Massey-Lie products \cite{HS1} which are
homotopy invariant.

Unfortunately, the Massey-Lie
products do not allow us to formulate a
criterion of formality of a DG algebraic structure.
Recently, D. Kaledin \cite{Kaledin} found
a good refinement of the notion of
Massey-Lie products which allowed him to
formulate the criterion of formality for
DG (associative) algebras.

While this manuscript was in preparation
paper \cite{Gilles} by G. Halbout appeared
on arxiv.org.
In this paper the author uses Tamarkin's technique
 \cite{Hi}, \cite{Dima} to establish
the formality of the DGLA of Drinfeld over any field of
characteristic zero. In \cite{Gilles} G. Halbout also used this
result to get a partial answer to Drinfeld's question \cite{D-Op}
about the quantization of coboundary Lie bialgebras. It would be
interesting to further develop this technique and get a complete
answer to this question.


\begin{thebibliography}{99}

\bibitem{Peter}  D. Ananikian, P.P. Kulish and V.D. Lyakhovsky,
Full chains of twists for symplectic algebras,
St. Petersburg Math. J. {\bf 14}, 3 (2003) 385--404; math.QA/0109083.

\bibitem{C} D. Calaque, Formality for Lie algebroids, Comm. Math. Phys. {\bf 257}, 3 (2005) 563--578;
math.QA/0404265.

\bibitem{Peter1} E. Celeghini, P. P. Kulish,
Deformation of orthosymplectic Lie superalgebra $osp(1 \vert 4)$, J.
Phys. A. {\bf 37}, 20  (2004) 211--216.

\bibitem{D} V.A. Dolgushev, Covariant and equivariant formality
theorems, Adv. Math.  {\bf 191}, 1  (2005) 147--177; math.QA/0307212.

\bibitem{D1} V.A. Dolgushev, A formality theorem for Hochschild
chains, math.QA/0402248.

\bibitem{D2} V.A. Dolgushev, A proof of Tsygan's formality conjecture for
an arbitrary smooth manifold, Ph. D. Thesis, MIT,   math.QA/0504420.

\bibitem{D83}  V.G. Drinfeld, Soviet. Math. Doklady {\bf 28}
(1983) 667-671.

\bibitem{D-Op} V.G. Drinfeld, On some unsolved problems in quantum group
theory, Lecture Notes in Math., {\bf 1510}, Springer, Berlin, 1992.

\bibitem{GJ} E. Getzler and J.D.S. Jones, Operads, homotopy algebra and iterated integrals for double loop spaces,
math.QA/9403055

\bibitem{GK} V. Ginzburg and M. Kapranov, Koszul duality for operads, Duke Math. J.  {\bf 76}, 1 (1994) 203-272.

\bibitem{Gilles} G. Halbout, Formality theorem for Lie bialgebras and quantization of coboundary $r$-matrices, math.QA/0506487.

\bibitem{HS1} S. Halperin and J. Stasheff, Obstructions to homotopy equivalences, Adv. Math.
{\bf 32}, 3 (1979) 233-279.

\bibitem{Hi} V. Hinich, Tamarkin's proof of Kontsevich formality theorem, Forum Math. {\bf 15},
4 (2003) 591-614; math. QA/0003052.

\bibitem{Kaledin} D. Kaledin, Some remarks on formality in families, math.AG/0509699.

\bibitem{Kontsevich} M. Kontsevich, Deformation quantization of
Poisson manifolds, Lett. Math. Phys. {\bf 66}, 3 (2003) 157-216.


\bibitem{Peter3}  P.P. Kulish, V.D. Lyakhovsky, and
M.A. del Olmo, Chains of twists for classical Lie algebras,
J. Phys. A: Math. Gen. {\bf 32}, 49 (1999) 8671 -- 8684;
math.QA/9908061.

\bibitem{Peter2}  P.P. Kulish, V.D. Lyakhovsky, and A. Stolin, Chains of
Frobenius subalgebras of so(M) and the corresponding twists,
J. Math. Phys.  {\bf 42}, 10  (2001)  5006--5019; math.QA/0010147

\bibitem{Q} D. Quillen, Homotopical algebra, Lecture Notes in Math. {\bf 43}
(1967).

\bibitem{Sh} B. Shoikhet, An explicit formula for the deformation
quantization of Lie bialgebras, math.QA/0402046.

\bibitem{Dima} D. Tamarkin, Another proof of M. Kontsevich formality
theorem, math.QA/9803025.

\bibitem{TT} D. Tamarkin and B. Tsygan, Cyclic formality and index theorems, Lett. Math. Phys. {\bf 56}
(2001) 85-97.
\end{thebibliography}
\end{document}